# A Graph-Theoretic Approach to Ring Analysis: Dominant Metric Dimensions in Zero-Divisor Graphs

**Nasir Ali[a], Hafiz Muhammad Afzal Siddiqui[a,] Muhammad Imran Qureshi[b]**

**[a]Department of Mathematics, COMSATS University Islamabad, Lahore Campus, Pakistan.**

**[b]Department of Mathematics, COMSATS University Islamabad, Vehari Campus, Pakistan.**

**Email address:**

Corresponding Author : nasirzawar@gmail.com (Nasir Ali), hmasiddiqui@gmail.com (Hafiz Muhammad Afzal Siddiqui), imranqureshi18@gmail.com (Muhammad Imran Qureshi),

**ORCID ID: https://orcid.org/0000-0003-4116-9673 (Nasir Ali),**

**Abstract:** This article investigates the concept of dominant metric dimensions in zero divisor graphs (ZD-graphs) associated with rings. Consider a finite commutative ring with unity, denoted as $R$, where nonzero elements $x$ and $y$ are identified as zero divisors if their product results in zero $(x.y = 0)$. The set of zero divisors in ring $R$ is referred to as $L(R)$. To analyze various algebraic properties of $R$, a graph known as the zero-divisor graph is constructed using $L(R)$. This manuscript establishes specific general bounds for the dominant metric dimension (Ddim) concerning the ZD-graph of $R$. To achieve this objective, we examine the zero divisor graphs for specific rings, such as the ring of Gaussian integers modulo $m$, denoted as $Z_m[i]$, the ring of integers modulo $n$, denoted as $\mathbb{Z}_n$, and some quotient polynomial rings. Additionally, we present a general result outlining bounds for the dominant metric dimension expressed in terms of the maximum degree, girth, clique number, and diameter of the associated ZD-graphs. Finally, we provide insights into commutative rings that share identical metric dimensions and dominant metric dimensions. This exploration contributes to a deeper understanding of the structural characteristics of ZD-graphs and their implications for the algebraic properties of commutative rings.

## 1. Introduction

Beck [1] proposed the connection between graph theory and algebra by introducing a ZD-graph of a commutative ring $R$. His focus was on node coloring in a graph, specifically exploring how ring elements correspond to these nodes. It's worth noting that in this context, a zero vertex is connected to all other vertices. The $Z^o(R)$ notation is commonly used to refer to this type of ZD-graph in literature.

In [2], Anderson and Livingston delved into a ZD-graph where each node represents a nonzero zero divisor (ZD). They defined an undirected graph where nodes $x$ and $y$ are connected by an edge if and only if $xy = 0$, naming it a ZD-graph of $R$. Anderson and Livingston primarily concentrated on finite rings, producing finite graphs when $R$ is finite. Their objective was to ascertain whether a graph is complete or a star for a given ring. Here, $\Gamma(R)$ is used to denote this type of ZD-graph of $R$. It's worth noting that this ZD-graph definition slightly differs from Beck's, as zero is not considered a vertex here. Consequently, $\Gamma(R) \subseteq Z^o(R)$. Anderson and Livingston [2] established a connection between the properties of $R$ and the graph properties of $\Gamma(R)$, yielding commutative ringucial insights into $\Gamma(R)$.

Expanding the ZD-graph concept from unital commutative rings to noncommutative rings, Redmond [3] introduced various methods to characterize ZD-graphs related to noncommutative rings. His approach encompassed both undirected and directed graphs. In a subsequent work [4], Redmond extended this concept to commutative rings, transforming it into an ideal-based ZD graph. The goal was to generalize the method by replacing elements with zero products with elements whose product belongs to a specific ideal $I$ of ring $R$.

An ideal-based ZD-graph denoted by $\Gamma_I(R)$ can be obtained by considering two nonzero ZD, $x$ and $y$ of $R$ as nodes, and there is an edge between them if and only if $xy \in I$ such that $\{xy \in I, \text{for some}, y \in R - I\}$. Various types of graphs, like total graphs, unit graphs, and Jacobson graphs, have been defined by different authors [5,6,7,8,9]. For ring theory basics, refer to [12, 13]. For graph theory basics, see [10, 11]. [24-25] discuss unique parameters for various graphs.

The graph linked with $R$ displays the features of $L(R)$. It visually and analytically uncovers the algebraic properties of rings using graph theory. In [5], authors explored ZD-graph properties. To study $\Gamma(R)$, we follow Anderson and Livingston's method [5], where non-zero zero divisors commutative rings are vertices of $\Gamma(R)$. In this paper, unless stated otherwise, we assume $R$ is a finite unital commutative ring. $L(R)$ is the set of non-zero zero divisors. If $R$ has only one maximal ideal, it's termed local.

For any $x \in R$, the annihilator of $x$, denoted as $ann(x)$, consists of all $y \in R$ such that $xy = 0$. An element $r$ in a ring $R$ is considered nilpotent if $r^m = 0$. A reduced ring is one that contains no nilpotent elements except for zero. The ring of integers modulo $n$, denoted as $\mathbb{Z}_n$, is the set $\{0, 1, 2, \ldots, n - 1\}$. Furthermore, $Z_n[i]$ denotes the ring of Gaussian integers modulo $n$ and, defined as $\{x + iy : x, y \text{ are in } \mathbb{Z}_n \text{ and } i^2 = -1\}$, under complex multiplication and addition. Here, $F$ signifies a finite field. Osba et al. introduced the graph for $Z_n[i]$ in [14]. The ZD-graph of the Gaussian integer ring $Z_n[i]$ is denoted as $\Gamma(Z_n[i])$.

The subset $S \subseteq V(G)$ is said to be a dominating set of $G$ if, for every vertex $x$ in $V(G) \setminus S$, there exists minimum one vertex $u \in S$ such that $x$ is adjacent to $u$. The set with minimum cardinality among the dominating set of $G$ is called a dominating number of $G$ and is denoted by $\gamma(G)$ [26]. An ordered set $W = \{w_1, w_2, \ldots, w_k\} \subseteq V(G)$ is called a resolving set of $G$ if every pair of vertices $u, v \in V(G)$ have distinct representation with respect to $W$, that is, $r(u|W) \neq r(v|W)$, where $r(u|W) = (d(u, w_1), d(u, w_2), \ldots, d(u, w_k))$. The minimum cardinality among the resolving set of $G$ is called the metric dimension (MD) of $G$ and is denoted by $dim(G)$. Brigham et al. [27] combined the concept of MD and dominating set by term resolving domination number, denoted by $\gamma_r(G)$ and got some result that $max\{dim(G), \gamma(G)\} \leq \gamma_r(G) \leq dim(G) + \gamma(G)$. Later, Henning and Oellarmann [28] studied the metric locating dominating number of graph $G$, denoted by $\gamma_M(G)$ and found lower and upper bounds, that is, $\gamma(G) \leq \gamma_M(G) \leq n - 1$. Then, Gonzalez et al. [29] examined different lower and upper bounds, i.e., $max\{dim(G), \gamma(G)\} \leq \gamma_M(G) \leq dim(G) + \gamma(G)$. An ordered set $W \subseteq V(G)$ is called a dominant resolving set of $G$ if $W$ is a resolving set and a dominating set of $G$. The dominant resolving set having minimum cardinality is called a dominant basis of $G$. In contrast, the cardinality of the dominant basis is called a dominant metric dimension of $G$ and is denoted by $Dim_d(G)$.

The motivation behind this study stems from a profound curiosity about the structural intricacies of zero divisor graphs associated with finite commutative rings. These graphs, derived from the concept of zero divisors in rings, offer a unique lens through which we can explore and understand the underlying algebraic properties of these mathematical structures. The focus on dominant metric dimensions adds a layer of complexity to our investigation, as it plays a pivotal role in delineating the metric characteristics of ZD-graphs.

We seek to establish specific bounds for the dominant metric dimension within the context of ZD-graphs. By doing so, we aim to contribute not only to the theoretical understanding of these graphs but also to their practical implications in various algebraic scenarios. The chosen rings for our analysis, including the ring of Gaussian integers $Z_m[i]$ modulo $m$, the ring of integers $\mathbb{Z}_n$ modulo $n$, and other quotient polynomial rings, serve as diverse examples to illustrate the versatility and applicability of our findings.

Furthermore, our study extends beyond individual rings to provide general bounds for the dominant metric dimension, intertwining it with fundamental parameters such as the maximum degree, girth, clique number, and diameter of ZD-graphs. This comprehensive approach enables us to uncover overarching principles that govern the metric characteristics of these graphs.

Ultimately, our exploration aims to enhance our understanding of the intricate relationships between ZD-graphs and the algebraic properties of commutative rings. By unraveling these connections, we hope to contribute valuable insights to the broader mathematical community and inspire further research in this fascinating field.

The article makes several notable contributions. Firstly, it extends the concept of dominant metric dimension (Ddim) to ZD-graphs, particularly focusing on ZD-graphs representing the ring $\mathbb{Z}_n$ of integers modulo $n$. This generalization encompasses diverse characterizations of rings based on the vertices in the ZD-graph, including metric and dominant metric representations. The research establishes that the Ddim serves as a robust tool for effectively characterizing rings based on their unique graph structures. Secondly, the article delves into the characterization of various rings through their ZD-graphs, demonstrating instances where the Ddim can be effectively bounded by the graph's diameter. Lastly, the article introduces a straightforward procedure for calculating the dominant metric dimension of ZD-graphs that represent rings of integers $\mathbb{Z}_n$ modulo $n$. The novelty of determining the Ddim of graphs lies in its ability to offer a more comprehensive understanding of both the structural and algebraic properties of graphs. Such insights hold significant relevance in practical applications, including network design, social networking, and communication systems, where a nuanced understanding of graph properties is commutative ringucial for optimal system design.

## 1.1 Preliminaries

Formally, the graph is an ordered pair $G = (V, E)$; here, $V$ and $E$ stands for vertices or nodes and edge set, respectively. A graph's order and size are defined as the cardinality of nodes set and edges set, respectively. The open neighborhood of a node $v$ is written as $N(v)$, and defined as $\{v \in V(G): vu \in E(G)\}$, while the closed neighborhood of a node $u$ is written as $N[u]$, and defined as $\{u\} \cup N(u)$. The distance $d(u', v')$ between two nodes $u'$ and $v'$ is defined as the length of the shortest path between them, while $d(w, e') = min\{d(w, u'), d(w, v')\}$ defines the distance between a node $w$ and the edge $e' = u'v'$.

The length of the longest path is the diameter of the graph, which is denoted by $diam(G)$. Mathematically, $diam(G) = sup\{d(r, s)$: where r and s are distinct vertices in G$\}$. Let $H$ be a subset of set of vertices along with any subset of edges containing those vertices is said to be a subgraph of a graph $G$; mathematically, it is denoted by $H \subset G$. Consider a smallest cycle subgraph $H$ in a graph $G$, then the number of edges in $H$ is called the girth of the graph, denoted by $gr(G)$. A clique is defined as maximal complete subgraph of a graph $G$ which is denoted by $K$ and $|K| = \omega(G)$ is called the clique number.

A graph is said to be a regular graph if for every $r \in V, deg(r) = c$ for a fix, $c \in Z^+$. A graph is said to be complete if there is a connection between all pairs of vertices. It is represented by $k_m$, where $m$ is the number of vertices. A graph is considered a complete bipartite graph if it can be divided into two distinct sets of vertices, $X$ and $Y$, where each vertex in $X$ is connected to every vertex in $Y$, and it is usually denoted by $k_{m,n}$, where $|X| = m$ and $|Y| = n$. When a vertex vanishes from a connected graph and creates two or more disconnected components of the graph, it is called a cut vertex.

Kelenc et al. [16] studied the edge metric-dimension (EMD) of various graphs, including the complete graph $k_m$, the path graph $P_n$, and complete bipartite graph $k_{m,n}$. The relationship between the MD and the EMD allows for the identification of graphs where these two dimensions are equal, as well as for some other graphs $G$ for which $dim(G) < dim_E(G)\ or\ dim_E(G) < dim(G)$. Basically, Kelenc et al. explored the comparison of values $dim(G)$ and $dim_E(G)$.

Recently, a study on metric parameters for ZD-graphs has been done. Redmond, in 2002 [17], studied the ZD-graphs of noncommutative rings, and in 2003 [18], the ideal-based ZD-graphs of commutative rings were studied by him. In 2019 [19], the metric dimension of ZD-graphs for ring $\mathbb{Z}_n$ was calculated. In 2020 [20], bounds for the EMD of ZD-graphs related to rings were studied by Siddiqui et al. Pirzada and Aijaz in 2020 [21] studied ZD-graphs for commutative rings for their metric and upper dimension.

Susilowati et al. [15] conducted a study on the dominant metric dimension of a specific class of graphs, providing valuable insights into the characterization of graphs with distinct dominant metric dimensions. Their research extended to determining the dominant metric dimension of joint and comb products of graphs. In this context, we build upon and consider pertinent results presented in [30-31]. These findings contribute to the ongoing exploration of dominant metric dimensions in graph theory, further enriching our understanding of the structural properties and characteristics of graphs.

1. For path graph denoted by $P_n$ and cyclic graph $C_n$, $\gamma(P_n) = \gamma(C_n) = \left\lceil \frac{n}{3} \right\rceil$ and $dim(P_n) = 1$ and $dim(C_n) = 2$.
2. For a complete graph denoted by $K_n$, $\gamma(K_n) = 1$ and $dim(K_n) = n - 1$.
3. For a start graph denoted by $S_n$, $\gamma(S_n) = 1$ and $dim(S_n) = n - 2, \forall\ n \geq 2$.
4. For a complete bipartite graph $K_{m,n}$, $\gamma(K_{m,n}) = 2$ and $dim(K_{m,n}) = m + n - 2,, \forall\ m,n \geq 2$.

In addition, we delve into prior findings on the dominant metric dimension of some graphs as presented in [15]. This exploration aims to build upon and incorporate insights gained from the research conducted by Susilowati et al. The consideration of these earlier results contributes to a more comprehensive understanding of the dominant metric dimension of graph $G$, adding depth to the existing body of knowledge in this:

**Theorem 1** [15]**:** If $C_n$ is a cyclic graph of order $n \geq 7$, then $Dim_d(C_n) = \gamma(C_n)$.
**Theorem 2** [15]**:** If $G$ is a start graph $S_n$, having order $n \geq 2$, then $Dim_d(G) = n - 1$.
**Theorem 3** [15]**:** Let $K_{m,n}$ be a complete bipartite graph with conditions that $m, n \geq 2$, then $Dim_d(K_{m,n}) = dim(K_{m,n})$.
**Theorem 4** [15]**:** If $G$ is a path graph $P_n$, with $n \geq 4$, then $Dim_d(P_n) = \gamma(P_n)$.
**Theorem 5** [15]**:** If $G$ is a complete graph $K_n$, with $n \geq 2$, then $Dim_d(K_n) = dim(K_n)$.

**Theorem 6** [15]:  $Dim_d(P_n) = 1 \Leftrightarrow G \cong P_n, n = 1,2$.

Additionally, it is noteworthy that the dominant metric dimension for a single-vertex graph G is considered to be zero. On the other hand, for an empty graph, the Ddim is undefined. With these considerations in mind, our discussion initiates with the following observation.

## 2.  Results

### 2.1.1.  *Dominant Metric dimension of some zero-divisor graphs*

**Theorem 2.1.** Consider a finite commutative ring  $R$  with unity. Then

i.  $Dim_d(\Gamma(R))$  is finite if and only if  $R$  is finite.

ii.  $Dim_d(\Gamma(R))$  is undefined if and only if  $R$  is an integral domain.

**Proof.**

i.  Suppose that  $Dim_d(\Gamma(R))$  is finite. In this case, one can find a minimal dominant metric basis for Γ(R), denoted as  $\{v_1, v_2, \ldots, v_t\}$ . Using ([2], Theorem 2.3), we establish that the  $diam(\Gamma(R))$  is bounded by 3, i.e.,  $diam(\Gamma(R)) \leq 3$ . Consequently, for every pair  $(r, e)$  where  $r$  belongs to the vertices  $V(\Gamma(R))$  and  $e$  belongs to the edges  $E(\Gamma(R))$ , the metric distance  $d(r, e)$  is limited to  $0, 1, 2, or\ 3$ . This restriction implies that the size of the line graph  $|L(R)|$  is at most  $4t$ . As a result, Γ(R) is finite, leading to the conclusion that  $R$  is also finite. Fore the converse part, given that  $R$  is finite, hence  $|L(R)|$  is also finite. Moreover, as Γ(R) is contained in  $R$ . It follows that,  $Dim_d(\Gamma(R))$  is finite.

ii.  As we know, if  $R$  is an integral domain, then Γ(R) is not defined, which shows that.  $Dim_d(\Gamma(R))$  is undefined and vice versa. ∎

The following outcome presented discloses the Ddim of the ZD-graph of ring  $R$  whenever  $\Gamma(R)$  is isomorphic to  $P_m$  for some integer  $m$ .

**Proposition 2.1.** Consider a finite commutative ring  $R$  with unity. Subsequently,  $Dim_d(\Gamma(R)) = 1$  if and only if  $R$  exhibits an isomorphism with any of the following rings:  $\mathbb{Z}_6, \mathbb{Z}_8, \mathbb{Z}_9, \mathbb{Z}_2 \times \mathbb{Z}_2, \mathbb{Z}_3(r)/(r^2), \mathbb{Z}_2(r)/(r^3), or\ \mathbb{Z}_4(r)/(2r, r^2 - 2)$ .

**Proof.** Assume that  $Dim_d(\Gamma(R)) = 1$ . In such cases, the ZD-graphs of the rings  $\mathbb{Z}_6, \mathbb{Z}_8, \mathbb{Z}_9, \mathbb{Z}_2 \times \mathbb{Z}_2, \mathbb{Z}_3(r)/(r^2), \mathbb{Z}_2(r)/(r^3), or\ \mathbb{Z}_4(r)/(2r, r^2 - 2)$  are path graphs.It is well-established by Theorem 6 that path graphs are the only graphs with $Dim_d = 1$. Furthermore, based on ([19], Lemma 2.6), if $\Gamma(R)$ is isomorphic to $P_n$ then $|L(R)|$ is maximum three.

Case I: If $\Gamma(R) \cong P_2$, then $|L(R)| = a, b$ such that $a.b = 0$. Rings that satisfy this property include $\mathbb{Z}_9, \mathbb{Z}_2 \times \mathbb{Z}_2, \mathbb{Z}_3(r)/(r^2)$.

Case II: If $\Gamma(R) \cong P_3$, then $|L(R)| = a, b, c$, such that $a.b = 0$ and $b.c = 0$. Then, rings that satisfy this property are $\mathbb{Z}_6, \mathbb{Z}_8, \mathbb{Z}_2(r)/(r^3)$, and $\mathbb{Z}_4(r)/(2r, r^2 - 2)$[14].

Conversely, the ZD-graphs of the rings mentioned above are isomorphic to either. $P_2 \; or \; P_3$, [14]. Hence, by Theorem **6**, $Dim_d(\Gamma(R)) = 1$. ∎

**Proposition 2.2.** Consider a finite commutative ring $R$ with unity. Let $R$ exhibits an isomorphism with any of the following rings, $\mathbb{Z}_3 \times \mathbb{Z}_3$, $\mathbb{Z}_2[r, s]/(r, s)^2$, $K_4(r)/(r^2), \mathbb{Z}_4(r)/(r^2 + r + 1), \mathbb{Z}_4(r)/(2, r)^2$. Then $\Gamma(R) \cong C_m$, and $Dim_d(\Gamma(R))$ is 2.

**Proof.** For given ring $R$, $\Gamma(R)$ is a cyclic graph, according to ([3], Theorem 2.4). Moreover, the length of the cyclic graph is at most 4. As $\Gamma(R)$ is isomorphic to $C_m$ and $m$ does not exceed 4. It follows that $Dim_d(\Gamma(R)) = 2$.

The corresponding ZD-graphs for the rings mentioned above can be observed in Figure 1.∎

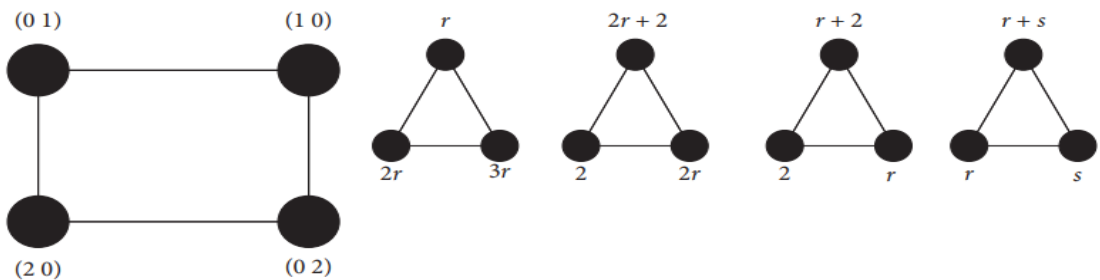

*Figure 1: ZD-graphs of $\mathbb{Z}_3 \times \mathbb{Z}_3$, $K_4(r)/(r^2), \mathbb{Z}_4(r)/(r^2 + r + 1), \mathbb{Z}_4(r)/(2, r)^2, \mathbb{Z}_2[r, s]/(r, s)^2$*

**Theorem 2.2.** Consider a finite commutative ring $R$ with unity. Moreover, if each $r \in L(R)$ is nilpotent then,

a) If $|L(R)| \geq 3$ and $L(R)^2 = \{0\}$, then $Dim_d(\Gamma(R)) = |L(R)| - 1$.

b) If $|L(R)| \geq 3$ and $L(R)^2 \neq 0$, then $Dim_d(\Gamma(R))$ is finite.

**Proof.**

a)  Given that $|L(R)| \geq 3$ and $L(R)^2 = \{0\}$ then $r.s = 0$ for all $r, s \in L(R)$ and by ([5], Theorem 2.8**)**, $\Gamma(R)$ is complete graph. Therefore, by applying Theorem 5 and Remark 2, $Dim_d(\Gamma(R)) = |L(R)| - 1$.

b)  Given that $|L(R)| \geq 3$ and $L(R)^2 \neq 0$, then it implies the existence of some $r \in L(R)$ such that $r^2 = 0$. This further implies that there exists $s \in L(R)$ such that $d(r, s) \geq 2$. Consequently, $L(R)/(r, s)$ serves as the dominant metric generator for any vertex $s$ adjacent to $r$. As a consequence, $Dim_d(\Gamma(R))$ is finite. ∎

**Theorem 2.3.** Consider a finite commutative ring $R$ with unity. Moreover, let $|L(R)| \geq 3$. Let the associated ZD-graph $\Gamma(R)$ has a cut vertex but have no vertex with degree 1, then $Dim_d(\Gamma(R))$ is 3 or 5.

**Proof.** For the given ring $R$. Let the associated ZD-graph $\Gamma(R)$ has a cut vertex but have no vertex with degree 1 then according to ([25], Theorem 3), $R$ exhibits an isomorphism with any of the following rings: $\frac{\mathbb{Z}_2[r,s]}{(r^2, s^2 - rs)}, \frac{\mathbb{Z}_4[r]}{(r^2 + 2r)}, \frac{\mathbb{Z}_4[r,s]}{(r^2, s^2 - rs, \ rs - 2, \ 2r, \ 2s)}, \frac{\mathbb{Z}_8[r,s]}{(2r, \ r^2 + 4)}, \frac{\mathbb{Z}_2[r,s]}{(r^2, s^2)}, \frac{\mathbb{Z}_4[r]}{(r^2)}, \frac{\mathbb{Z}_4[r,s]}{(r^2, s^2, \ rs - 2, 2r, 2s)}$.

Figure 2(a) which represents the $\Gamma(R)$ for the first four rings, the set $D_1 = \{v_1, v_3, v_4, v_5, v_7\}$ is identified as a minimum dominant metric generator for $\Gamma(R)$.

In Figure 2(b), representing $\Gamma(R)$ for the remaining three rings, the set $D_2 = \{x_1, x_2, x_4\}$ is recognized as a minimum dominant metric generator for $\Gamma(R)$. Therefore, it is concluded that $Dim_d(\Gamma(R))$ can either be 3 or 5 based on the minimum dominant metric generators identified for the respective figures. ∎

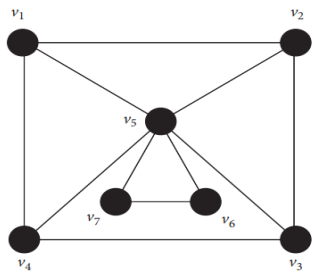 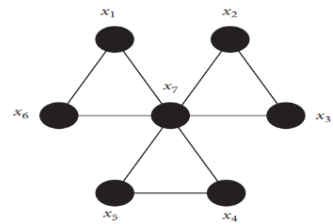

*Figure 2 (a)*                    *Figure 2 (b)*

**Theorem 2.4.** Let $R$ be a finite commutative ring with unity and $R \cong \mathbb{Z}_2 \times \mathbb{K}$, here, $\mathbb{K}$ denotes some field. Then $\Gamma(\mathrm{R}) \cong K_{1,|L(R)|-1}$ and $Dim_d\big(\Gamma(\mathrm{R})\big) = |L(R)| - 1$. Furthermore, Let $R$ be a local ring having no cycles in associated $\Gamma(\mathrm{R})$ then $Dim_d\big(\Gamma(\mathrm{R})\big) = 1$.

**Proof.** Consider $R$ is a non-local ring such that $R \cong \mathbb{Z}_2 \times \mathbb{K}$. Then $L(R) = \{(1,0), (0,r) : r \in \{1,2,\ldots,|\mathbb{K}|-1\}\}$ such that $(1,0).(0,r) = 0, \forall\, r \in \{1,2,\ldots,|\mathbb{K}|-1\}$. Note that the vertex $(1,0)$ is the central vertex, which is adjacent to all other $|L(R)| - 1$ vertices. Hence, $\Gamma(\mathrm{R}) \cong K_{1,|L(R)|-1}$, so by Theorem 2, $Dim_d\big(\Gamma(\mathrm{R})\big) = |L(R)| - 1$. Further, consider $R$ is a local ring having no cycle in associated $\Gamma(\mathrm{R})$, then by ([25], *Theorem* 2.1), $\Gamma(\mathrm{R})$ is isomorphic to $P_2$ or $P_3$. Hence $Dim_d\big(\Gamma(\mathrm{R})\big) = 1$. ∎

**Theorem 2.6.** Consider the ring of integers $\mathbb{Z}_n$ modulo $n$. Assuming $p$ and $q$ as distinct primes, we have $Dim_d(\Gamma(\mathbb{Z}_n))$ as follows:

| $n$ | $Dim_d(\Gamma(\mathbb{Z}_n))$ |
|:---:|:---:|
| $p$ | $undefined$ |
| $p^2$ | $p - 2$ |
| $pq$ | $p + q - 4$ |
| $2^2$ | $0$ |
| $2^3$ | $1$ |
| $3^2$ | $1$ |

**Proof.** To establish the validity of the theorem, it is essential to meticulously examine and address each case separately.

(i)   When $n = p$, then $\Gamma(\mathbb{Z}_n)$ is an empty graph. So $Dim_d(\Gamma(\mathbb{Z}_n))$ is undefined.

(ii)   When $n = p^2$, then $\Gamma(\mathbb{Z}_n)$ consists of the complete graph $K_{p-1}$. By Theorem 5 and Remark 2, $Dim_d(\Gamma(\mathbb{Z}_n)) = p - 2$

(iii)  Now, if $n = pq$, where $p$ and $q$ are distinct primes, then we can partition the vertices into sets $U = \{\lambda p \in Z(\mathbb{Z}_n); (\lambda, q) = 1\}$ and $V = \{\lambda q \in Z(\mathbb{Z}_n); (\lambda, p) = 1\}$, clearly, partition shows that $\Gamma(\mathbb{Z}_n)$ is bipartite; moreover, we see $x.y = 0$ for all $x \in U$ and $y \in V$. Hence $\Gamma(\mathbb{Z}_n)$ is a complete bipartite graph. So, we conclude that when $n = pq$, then $\Gamma(\mathbb{Z}_n) \cong K_{q-1,p-1}$. Hence, by Theorem 3 and Remark 4, $Dim_d(\Gamma(\mathbb{Z}_n)) = p + q - 4$.

(iv)  When $n = 2^2$, then $\Gamma(\mathbb{Z}_n)$ is a single vertex graph. So $Dim_d(\Gamma(\mathbb{Z}_n))$ is zero.

(v)  When $n = 2^3$, then $\Gamma(\mathbb{Z}_n)$ is a path graph with two vertices. So, by Theorem 6, $Dim_d(\Gamma(\mathbb{Z}_n))$ is 1.

(vi)  When $n = 3^2$, then $\Gamma(\mathbb{Z}_n)$ is a path graph with two vertices. So, by Theorem 6, $Dim_d(\Gamma(\mathbb{Z}_n))$ is 1. Which completes the proof. ∎

One can see Table 1 below to understand the possible structures of ZD-graphs of $\mathbb{Z}_n$ for different values of $n$.

**Table 1**: *Dominant Metric dimension of $\Gamma(\mathbb{Z}_n)$*

| N | \|V\| | \|E\| | Diameter | Girth | $\Gamma(\mathbb{Z}_n)$ | $Dim_d(\Gamma(\mathbb{Z}_n))$ |
|---|---|---|---|---|---|---|
| $P$ | 0 | 0 | 0 | Undefined | $\Gamma(\mathbb{Z}_n) = \emptyset$ | *undefined* |
| $2^2$ | 1 | 0 | 0 | Undefined | ● | 0 |
| $3^2$ | 2 | 1 | 1 | Undefined | ●——● | 1 |
| $P^2$ $p \geq 5$ | p-1 | $\left(\dfrac{p-1}{2}\right)$ | 1 | 3 | Complete graph $K_{p-1}$ | $p-2$ |
| $2^3$ | 3 | 2 | 2 | Undefined | ●——●——● | 1 |
| $pq$ | q-1+p-1 | (q-1)(p-1) | 2 | 4 | Complete bipartite graph $K_{q-1,p-1}$ | $p+q-4$ |

### 2.1.2.  *Bounds between dominant metric dimension and diameter of Zero-divisor graphs:*

In this section, our attention is directed towards establishing bounds for the dominant metric dimension,

girth, and clique number of ZD-graphs.

**Theorem 2.1.2.1** Consider a finite commutative ring $R$ with unity. Let $R \cong \mathbb{K}_1 \times \mathbb{K}_2$, where $\mathbb{K}_1 \& \mathbb{K}_2$ denotes some finite fields with $|\mathbb{K}_1| = m \geq 3$, $|\mathbb{K}_2| = n \geq 3$. Then, $Dim_d\big(\Gamma(R)\big) = |\mathbb{K}_1| + |\mathbb{K}_2| - gr(\Gamma(R))$.

**Proof.** When $R \cong \mathbb{K}_1 \times \mathbb{K}_2$, then each vertex of the form $(u, 0) \in Z(R)$ is adjacent to each vertex of the form $(0, v)$ and Contrariwise. So $L(R)$ can be divided into two separate sets in the following manner: $\quad U = \{(u, 0), u \in \mathbb{K}_1\} \text{ and } V = \{(0, v), v \in \mathbb{K}_2\}$. Hence, $\quad \Gamma(R) \cong K_{m-1, n-1}$ and $gr\big(\Gamma(R)\big) = 4$, so by Theorem *3*, along with Remark 4, $Dim_d\big(\Gamma(R)\big) = |\mathbb{K}_1| + |\mathbb{K}_2| - gr(\Gamma(R))$. ∎

**Theorem 2.1.2.2** Consider a reduce ring $R$ and let $I_1$ and $I_2$ be two ideals such that $I_1 \cap I_2 = \{0\}$. Then $Dim_d\big(\Gamma(R)\big) = |I_1| + |I_2| - 2\,\omega\big(\Gamma(R)\big)$.

**Proof.** By using (*[20]*, Proposition *5*), $\Gamma(R) \cong K_{|I_1|, |I_2|}$. Also $\omega\big(\Gamma(R)\big) = 2$. Hence, by Theorem 3, along with Remark 4, $Dim_d\big(\Gamma(R)\big) = |I_1| + |I_2| - 2\,\omega\big(\Gamma(R)\big)$. ∎

**Lemma 2.1.2.1** Consider a ZD-graph for a ring $R$ having diameter at most 2 and $\Gamma(R) \ncong P_n$, then $Dim_d(\Gamma(R))$ is finite.

**Proof.** Consider the ZD-graph of diameter at most 2 and $\Gamma(R)$ is a non-path graph, then it means that $\Gamma(R)$ can be a cycle with at most 5 vertices, stars, a complete graph, or a Peterson graph. Then, using [Theorem 1-6], the result follows. ∎

Now, Let us determine the dominant metric dimension of the zero-divisor graph of the ring of Gaussian integers $Z_m[i]$ modulo $m$. As stated above, a ring of Gaussian integers modulo $m$ is denoted and defined as $Z_m[i] = \{x + iy : x, y \in Z_m \text{ and } i^2 = -1\}$ under the complex multiplication and addition, and $F$ denotes finite field. A Gaussian prime is a prime element in $Z[i]$, and the Gaussian primes can be characterized as follows:

(1) The elements $1 - i$ and $1 + i$ are Gaussian primes.

(2) Let $q$ be a prime integer such that $q \equiv 1\,(mod\ 4)$ and $q = a^2 + b^2$ for some

integers $a$ and $b$, then $a + ib$ and $a - ib$ are Gaussian primes.

(3) Let $p$ be a prime integer such that $p \equiv 3 \ (mod \ 4)$, then $p$ is a Gaussian prime. Additionally, for a Gaussian prime $q$, we have $-q, iq$, and $-iq$ Gaussian primes that is its complex conjugate and its associates are also Gaussian primes. Let $n = p \equiv 3 \ (mod \ 4)$ , then $Dim_d(\Gamma(Z_m[i]))$ is undefined since $Z_m[i]$ is a field, $\Gamma(Z_m[i])$ is an empty graph.

**Theorem 2.1.2.3.** Let $R$ be a finite commutative ring with unity and $R \cong Z_m[i]$. Then,

(1) For $n = p^2$, $Dim_d(\Gamma(Z_m[i])) = p^2 - 2$.

(2) For $p_j \equiv 3 \ (mod \ 4)$, $j = 1, 2$, then $Dim_d(\Gamma(Z_m[i])) = p_1{}^2 - p_2{}^2 - 2\omega. \Gamma(Z_{p_1,p_2}[i])$.

(3) For $n = p \equiv 1 \ (mod \ 4)$ with $p = a^2 + b^2$, $Dim_d(\Gamma(Z_m[i])) = 2p - gr(\Gamma(Z_m[i]))$.

**Proof.**

(1) From ([14], Theorem 15), it can be seen that $\Gamma(Z_m[i])$ is a complete graph when $n = p^2$. So, by Theorem 5, along with Remark 2, $Dim_d(\Gamma(Z_m[i])) = p^2 - 2$.

(2) Consider $p_1, p_2$ are primes such that $p_j \equiv 3 \ (mod \ 4)$, where $j = 1, 2$, then $\Gamma(Z_{p_1,p_2}[i])$ is isomorphic to a complete bipartite graph. As $Z_{p_1,p_2}[i] \cong Z_{p_1}[i] \times Z_{p_2}[i]$ Moreover, in the case of a complete bipartite graph, the $\omega$ is $2$. Therefore, according to Theorem 3, along with Remark 4, $Dim_d(\Gamma(Z_m[i])) = p_1{}^2 - p_2{}^2 - 2\omega. \Gamma(Z_{p_1,p_2}[i])$.

(3) If $n = p \equiv 1 \ (mod \ 4)$ with $p = a^2 + b^2$, then $\Gamma(Z_p[i]) \cong K_{p-1,p-1}$. The girth of a complete bipartite graph is $2$; hence, by Theorem 3, along with Remark 4, $Dim_d(\Gamma(Z_m[i])) = 2p - gr(\Gamma(Z_m[i]))$. ∎

Finally, we present a table having the rings with the same metric dimension and dominant metric dimensions.

*Table 2: Rings having same Metric Dimension and Dominant Metric dimension*

| Commutative Rings $R$ | $\dim(\Gamma(R)) = Dim_d(\Gamma(R))$ |
|---|---|
| $\mathbb{Z}_6, \mathbb{Z}_8, \mathbb{Z}_9, \mathbb{Z}_2 \times \mathbb{Z}_2, \mathbb{Z}_3(r)/(r^2), \mathbb{Z}_2(r)/(r^3), or \ \mathbb{Z}_4(r)$ $/(2r, r^2 - 2)$ | 1 |
| $\mathbb{Z}_3 \times \mathbb{Z}_3, \ K_4(r)/(r^2), \mathbb{Z}_4(r)/(r^2 + r + 1), \mathbb{Z}_4(r)/$ | 2 |
| $(2, r)^2, \mathbb{Z}_2[r, s]/(r, s)^2$ | |
| $\mathbb{K}_1 \times \mathbb{K}_2$, Where $\mathbb{K}_1$ and $\mathbb{K}_2$ are finite fields of order greater than or equal to $3$ | $|\mathbb{K}_1| + |\mathbb{K}_2| - 4$ |
| $\frac{\mathbb{Z}_p[\mathrm{r}]}{(r^2)}$ $or$ $\mathbb{Z}_{p^2}$, where $p$ is a prime | $p - 2$ |

## 3. Conclusions

In this work, the rings have been characterized by studying the dominant metric dimension of associated ZD-graphs. Dominant metric dimension of graphs linked to commutative rings, including the ring $\mathbb{Z}_n$ of integers modulo $n$, polynomial rings, and the ring of Gaussian integers modulo $m$ have been examined. Furthermore, the dominant metric dimension for the ring $\mathbb{Z}_n$ of integers, modulo $n$ is generalized. The study concluded by presenting bounds between Ddim, girth, clique number, and diameter of ZD-graphs. By delving into this research, we pave the way for researchers to explore a rich interdisciplinary landscape, offering promising insights that could have far-reaching impacts in diverse domains like network analysis and commutative ring cryptography.

**Declaration.**

- **Availability of data and materials:** The data is provided on request to the authors.

- **Conflicts of interest:** The authors declare that they have no conflicts of interest and all agree to publish this paper under academic ethics.

- **Fundings:** This received no specific grant from any funding agency in the public, commercial, or not-for-profit sectors.

- **Author's contribution:** All the authors equally contributed towards this work.